\documentclass{article}

\usepackage{amsmath} 
\usepackage{amssymb}
\usepackage{amsthm}
\usepackage{enumerate}

\usepackage{tikz}

\theoremstyle{definition}
\newtheorem{defi}{Def.}
\theoremstyle{plain}
\newtheorem{lemma}[defi]{Lemma}
\newtheorem{theorem}[defi]{Theorem}
\newtheorem{prop}[defi]{Proposition}
\newtheorem{kor}[defi]{Corollary}
\theoremstyle{remark}
\newtheorem{rem}[defi]{Remark}
\newtheorem{exa}[defi]{Example}

\newcommand{\R}[1]{\mathbb{R}^{#1}}
\newcommand{\prf}[1][]{\noindent {\bf Proof{#1}:} }

\DeclareMathOperator{\Cov}{Cov}

\DeclareMathOperator{\Var}{Var}
\DeclareMathOperator{\Lip}{Lip}

\title{A functional central limit theorem for integrals of stationary mixing random fields}
\author{J\"urgen Kampf and Evgeny Spodarev\\
{\small Institute of Stochastics, Ulm University, Helmholtzstr.\ 18, 89081 Ulm, Germany.}\\
 {\small juergen.kampf@uni-ulm.de}, {\small evgeny.spodarev@uni-ulm.de}
}
\begin{document}

\maketitle

\begin{abstract}
We prove a functional central limit theorem for integrals $\int_W f(X(t))\, dt$, where $(X(t))_{t\in\mathbb{R}^d}$ is a stationary mixing random field and the stochastic process is indexed by the function $f$, as the integration domain $W$ grows in Van Hove-sense.  We discuss properties of the covariance function of the asymptotic Gaussian process.
\end{abstract}


\section{Introduction}

A random field is a collection of random variables, indexed by the points of the Euclidean space. 
Random fields have applications in various branches of science, e.g.\ in medicine \cite{ATW, TW}, in geostatistics \cite{ChDe, Wack} or in materials science \cite{McSt, Tor}.  

Central limit theorems for Lebesgue integrals have been studied for a long time. In the 1970's first central limit theorems for integrals of the form $\int_{W_n} X(t) \, dt $ were shown \cite{BulZ,Leon}, where $(X(t))_{t\in\R{d}}$ is a random field and the integration domains $W_n$ tend to $\R{d}$ in an appropriate way (see Section \ref{s:mult}). 
Meschenmoser and Shashkin \cite{MS} showed a functional central limit theorem for Lebesgue measures of excursion sets of random fields, where the stochastic process is indexed by the level of the excursion set. The Lebesgue measure of the excursion set equals $\int_{W_n} \mathbf{1}_{[u,\infty)}(X(t))\, dt$, where $u$ is the level. We will extend the result of \cite{MS} to a functional central limit theorem for $\int_{W_n} f(X(t))\, dt$, where the stochastic process is indexed by the function $f$, which is assumed to be Lipschitz continuous. While replacing indicator functions by Lipschitz continuous functions is straight-forward, we need an entirely different approach, since the index set of the stochastic process is much larger now (it is the real line in \cite{MS} and the space of Lipschitz continuous functions in the present paper). 
For a survey on other limit theorems for random fields, see \cite{Spo}.

We will study the covariance function of the asymptotic Gaussian process. This covariance function is a symmetric non-negative definite bilinear form. For a certain class of random fields we will present infinite sequences of functions which are orthogonal w.r.t.\ this bilinear form. While for this class of functions we can show that the asymptotic variance vanishes only for functions that are constant a.e.\ (w.r.t.\ the marginal distribution of the random field), we will construct a non-trivial random field for which this bilinear form vanishes identically.   

This paper is organized as follows. In Section \ref{s:prel}, we collect prelimaries about mixing random fields and orthogonal polynomials, in particular L\'evy-Meixner systems. We prove a multivariate version of the announced central limit theorem in Section \ref{s:mult}. Section \ref{s:cov} is devoted to the examination of the covariance function of the limiting Gaussian process. Since all results derived in this section can be formulated in terms of the asymptotic covariance matrix of the multivariate central limit theorem, it is no problem to do this examination before the functional version is obtained. In Section \ref{s:func}, we derive the functional central limit theorem.

\section{Preliminaries}\label{s:prel}

In this section, we introduce tools from the theory of mixing random fields and from the theory of orthogonal polynomials which we will need in later sections.  

\subsection{Mixing concepts}

Let $(\Omega, \mathcal{A}, \mathbb{P})$ be a probability space and let $\mathcal{F}, \mathcal{G} \subseteq\mathcal{A}$ be two sub-$\sigma$-algebras. Then we define the $\alpha$-mixing coefficient by
\begin{align*} 
\alpha(\mathcal{F},\mathcal{G}) &:= \sup\{ |\mathbb{P}(F \cap G) -\mathbb{P}(F)\mathbb{P}(G)| \mid F\in\mathcal{F}, G \in\mathcal{G} \}. 
\end{align*} 

For an $\mathbb{R}^s$-valued random field $(X(t))_{t\in\R{d}}$  we put
\begin{align*}
\alpha_\gamma (r)=\sup\{ \alpha(& \sigma(X_I),\sigma(X_J)) \mid I, J \subseteq \R{d},\quad
\mbox{there is $u\in\R{d}$ with $\|u\|=1$ such that }\\ 
&\min\{\langle u,t\rangle\mid t \in I\} -\max\{\langle u,t\rangle\mid t \in J\}>r,\quad  V_d((I\cup J)+B^d) \le \gamma \},
&  r\ge 0, \gamma\ge 2\kappa_d. 
\end{align*}
Here $\sigma(X_I)$ is the $\sigma$-algebra generated by the random variables $X_t,\, t\in I$. The Minkowski sum $I+B^d$ is  defined by $I+B^d = \{ i+b \mid i\in I,\, b\in B^d\}$, $B^d:=\{x\in\R{d}\mid \|x\|\le 1\}$ is the $d$-dimensional closed Euclidean unit ball and $\kappa_d:=\lambda_d(B^d)$. In words, $I$ and $J$ lie in halfspaces separated by a strip of width $r$ and the parallel volume of $I$ and $J$ at distance $1$ does not exceed $\gamma$. 




\subsection{L\'evy-Meixner systems of orthogonal polynomials}\label{ss:Meix}

We consider a family of probability measures $\Psi_\lambda,\, \lambda\in(0,\infty),$ on $\mathbb{R}$ such that all moments of these probability measures exist. 
For each $\lambda\in(0,\infty)$ we let $x \mapsto Q_n(x;\lambda), \, n \in \mathbb{N}_0$, where $\mathbb{N}_0=\mathbb{N}\cup \{0\}$, 
denote the sequence of real-valued orthogonal polynomials w.r.t.\ $\Psi_\lambda$ such that $x\mapsto Q_n(x;\lambda)$ has degree $n$ and its leading coefficient is $1$. In particular, $Q_0(x,\lambda) \equiv 1$. Such a sequence exists if and only if the measure $\Psi_\lambda$ is not concentrated on finitely many points. Indeed, if $\Psi_\lambda$ is concentrated on $n$ points, $n\in\mathbb{N}$, then the restrictions of more than $n$ functions to these points cannot be linearly independent and hence more than $n$ functions cannot be orthogonal. On the other hand, no polynomial but the zero polynomial has norm $0$ if $\Psi_\lambda$ is not concentrated on finitely many points. Hence, in this case, the desired sequence can be obtained e.g.\ by applying the Gram-Schmidt-procedure to $1, x, x^2, \dots$. Obviously the polynomials $Q_n(x;\lambda)$ are determined uniquely. For a general introduction to orthogonal polynomials see \cite{Sz} or \cite{Sch}.  

\begin{defi}\label{d:Meixner} A system  $x\mapsto Q_n(x;\lambda)$, $\lambda\in(0,\infty), n\in\mathbb{N}_0,$ of polynomials is called \emph{L\'evy-Meixner system} if
 \begin{enumerate}[(i)] 
\item for each $\lambda\in (0,\infty)$ there is a probability measure $\Psi_\lambda$ on $\mathbb{R}$ such that $x\mapsto Q_n(x;\lambda),$ $n\in\mathbb{N}_0,$ are orthogonal w.r.t.\ $\Psi_\lambda$,
\item there are two open sets $U,V\subseteq \mathbb{R}$ containing $0$ and analytical functions $a:V\to U$, $b:U\to\mathbb{R},$ that fulfill $a(0)=0$, $a'(0)\ne 0$ and $b(0)=1$ such that
\[\sum_{n\in \mathbb{N}_0} Q_n(x;\lambda) \frac{z^n}{n!} = b^{-\lambda}(a(z)) \exp(x\cdot a(z)), \quad x\in\mathbb{R},\ \lambda>0, \ z\in V.\]
\end{enumerate} 
\end{defi}

The measures $\Psi_\lambda$ appearing in this definition are determined uniquely if their moment generating function exists on a neighborhood of $0$. Indeed, two probability measures on $\mathbb{R}$ which have the same system of orthogonal polynomials must also have the same sequence of moments, since the monomials $x^n, n\in\mathbb{N}_0,$ can be written as linear combinations of orthogonal polynomials in a unique way and the moments can be identified as appropriate coefficients of this linear combinations. It is well known (see e.g.\ \cite[\S\ II.5]{Fr}) that a probability measure is uniquely determined by its sequence of moments if its moment generating function is finite on a neighborhood of $0$.  

 Since the probability measures $\Psi_\lambda$ are determined uniquely, we may call the sets $\{\Psi_\lambda\mid \lambda\in(0,\infty)\}$ L\'evy-Meixner systems as well.

By \cite[p.\ 52]{Sch} we have the following lemma:
\begin{lemma}\label{l:Meixp}
Under the assumptions of Definition \ref{d:Meixner}, the moment generating function of $\Psi_\lambda$ is  $t\mapsto b^\lambda(t)$. 
\end{lemma} 

We conclude from this lemma that probability distributions corresponding to L\'evy-Meixner systems are always infinitely divisible. 

\begin{theorem}\label{T:CovMeix}
\item Let $Q_n(x;\lambda), \, n\in \mathbb{N}_0,\, \lambda\in (0,\infty),$ be a L\'evy-Meixner system. Then:
\begin{enumerate}[a)]
\item For all $n\in \mathbb{N}_0$ and $x,y\in\mathbb{R}$ and $\lambda_1,\lambda_2\in(0,\infty)$ we have
\[ Q_n(x+y;\lambda_1+\lambda_2) = \sum_{i=0}^n \binom{n}{i} Q_i(x;\lambda_1) Q_{n-i}(y;\lambda_2).\]
\item Let $V_i,i\in J$, be a finite collection of independent random variables $V_i\sim \Psi_{\lambda_i}$. For any subsets $A\subseteq J$ and $B\subseteq J$ we have
			\begin{itemize} 
			\item $\sum_{i\in A} V_i \sim \Psi_{\lambda_A},$ for $\lambda_A:=\sum_{i\in A} \lambda_i $, and
			\item $\Cov( Q_n(\sum_{i\in A} V_i; \lambda_A), Q_m(\sum_{i\in B} V_i; \lambda_B)) = 0$ for $m\ne n$.
			\end{itemize}
\end{enumerate}
\end{theorem}
See \cite{L75} and the literature cited therein.

The L\'evy-Meixner systems will turn out to belong either to one of four well-known families of probability distributions or to one exotic family. These exotic distributions are called \emph{Meixner cosine hyperbolic}-distribution, $\mathbf{Mch}(a,\mu)$, $a\in(-\pi,\pi),\, \mu>0,$ and their characteristic function is given by
\[ \phi(u)= \big(\frac{\cos(a/2)}{\cosh((u-ia)/2)}\big)^{2\mu}. \] 
Schoutens \cite[p.\ 57]{Sch} shows that these distributions have densities
\[ f(x) = \frac{ (2\cos(a/2))^{2\mu} }{2\pi\Gamma(2\mu)} \exp(ax) |\Gamma(\mu+ix)|^2,\quad x\in\mathbb{R}.\]


An affine transformation of a L\'evy-Meixner system is again a L\'evy-Meixner system. More precisely we have:
\begin{lemma}\label{L:affMeix}
Let $x\mapsto Q_n(x;\lambda), \, n\in \mathbb{N}_0, \lambda\in(0,\infty),$ be a L\'evy-Meixner system and let $m,c\in\mathbb{R}$, $m\ne 0,$ be constants. Then $x\mapsto Q_n(mx+c\lambda;\lambda), \,n\in \mathbb{N}_0, \lambda\in(0,\infty),$ is also a L\'evy-Meixner system.
\end{lemma}

\prf We have
\begin{align*}
\sum_{n=0}^\infty Q_n(mx+c\lambda;\lambda) \frac{z^n}{n!} &= \frac{1}{b(a(z))^\lambda} \exp\big( (mx+c\lambda) \cdot a(z)\big) = \frac{1}{\tilde b(\tilde a(z))^\lambda} \exp\big( x \tilde a(z)\big),
\end{align*}
where $\tilde a(z):=m a(z)$ and $\tilde b(t):= b(\frac{1}{m} \cdot t)\cdot \exp(-ct/m)$. \qed


The following theorem is shown in \cite[Sec.\ 4.2 and 4.3]{Sch}. 

\begin{theorem}\label{t:LMeixcl}
Up to the shifts and scalings described by Lemma \ref{L:affMeix} each L\'evy-Meixner system is one of the following:

{ \begin{tabular}{ll|l}
\multicolumn{2}{l|}{Distribution} & Polynomials \\ \hline
Normal & $\mathcal{N}(\mu,\sigma^2)$ (with constant ratio $\mu/\sigma$, & Hermite\\[-0.5 mm]
& parameterized by $\sigma$) \\[1mm]
Gamma & $\Gamma(1,\alpha)$ (with constant scale factor $\lambda=1$, & (generalized) \\[-0.5 mm]
& parameterized by $\alpha-1$)&Laguerre\\[1mm]
Poisson & $Pois(\lambda)$ & Charlier\\[1mm]
Pascal & $P(\gamma,\mu)$ (with constant single-probability-para-  & Meixner type-I\\[-0.5 mm]
&meter $\gamma$, parameterized by the size-parameter $\mu$) \\[1mm]
Meixner cosine &$Mch(a,\mu)$ (with constant $a$, parameterized by $\mu$) & Pollaczek\\[-0.5 mm]
 hyperbolic&&\\
\end{tabular}}
\end{theorem}

\section{The multivariate central limit theorem}\label{s:mult} 

In this section, we derive a multivariate central limit theorem for Lebesgue integrals of random fields.  


A sequence $(W_n)_{n\in\mathbb{N}}$ of compact subsets of $\mathbb{R}^d$ is called \emph{Van Hove-growing} (VH-growing) if
\[ \lim_{n\to\infty} \lambda_d(\partial W_n+B^d)/\lambda_d(W_n) =0, \]
where $\partial W$ denotes the boundary of $W$.

\begin{theorem}\label{T:multCLTmix}
Let $(X(t))_{t\in\mathbb{R}^d}$ be a stationary, measurable $\mathbb{R}$-valued random field. Let $f_1,\dots, f_s:\mathbb{R}\to\mathbb{R}$ be measurable functions. Assume that:
\begin{enumerate}[(i)]
\item There is some $\delta>0$ with $\mathbb{E} |f_i(X(0))|^{2+\delta}<\infty,\,i\in\{1,\dots, s\}.$
\item $ \int_{\R{d}} |\Cov( f_i(X(0)), f_j(X(t)))|\,dt<\infty,\quad i,j=1,\dots, s$. 
\item There are $n\in\mathbb{N}$ and $C, l>0$ with $n/d>l + (2+\delta)/\delta$ such that $\alpha_\gamma(r)\le C r^{-n}\gamma^l$ for all $\gamma\ge 2\kappa_d$ and $r>0$, where $a_\gamma(r)$ is the $\alpha$-mixing coefficient for $(X(t))_{t\in\mathbb{R}^d}$. 
\end{enumerate}
Let $(W_n)_{n\in\mathbb{N}}$ be a VH-growing sequence of compact subsets of $\R{d}$. 

Then
\[  \big(\Phi_n(f_1), \dots, \Phi_n(f_s)\Big) \stackrel{n\to\infty}{\longrightarrow} \mathcal{N}(0,\Sigma) \]
in distribution, where $\Sigma$ is the matrix with entries $\sigma_{ij}:=\int_{\R{d}} \Cov\big( f_i(X(0)), f_j(X(t))\big)\,dt,$ $i,j=1,\dots, s,$ and
\[ \Phi_n(f)= \frac{\int_{W_n} f(X(t)) \, dt - \lambda_d(W_n)\cdot \mathbb{E}\, f(X(0))}{\sqrt{\lambda_d(W_n)}}. \] 
\end{theorem}

\prf For $u=(u_1,\dots,u_d)\in\mathbb{R}^d$ we define an $\mathbb{R}$-valued random field $(Y_u(t))_{t\in\mathbb{R}^d}$ by
\[ Y_u(t):= \sum_{i=1}^s u_if_i(X(t)),\qquad t\in\mathbb{R}^d. \]
Now  \cite[Remark 1]{Gor} tells us 
\[\frac{ \int_{W_n} Y_u(t) \, dt - \lambda_d(W_n)\cdot \mathbb{E}\, Y_u(0)}{\sqrt{\lambda_d(W_n)}} \to \mathcal{N}(0,\sigma_u^2), \]
where
\[ \sigma_u^2:= \int_{\R{d}} \Cov\big(Y_u(0), Y_u(t)\big) \, dt = u^T \Sigma u.\]
So the Theorem of Cram\'er and Wold implies the assertion. \qed\medskip

\section{The asymptotic covariance}\label{s:cov}

In this section we examine the asymptotic covariance matrix $\Sigma$ appearing in Theorem \ref{T:multCLTmix}. 

\subsection{Diagonal form}\label{ss:Diag}

First we discuss how the functions $f_1,\dots, f_s$ have to be chosen (depending on the distribution of the random field $(X(t))_{t\in\mathbb{R}^d}$) such that the asymptotic covariance matrix $\Sigma$ becomes a diagonal matrix, i.e.\ that
\begin{equation} \int_{\R{d}} \Cov\big(f_i(X(0)), f_j(X(t)) \big) \, dt =0,\quad i\ne j.\label{e:Cov0} \end{equation}
 If $\Sigma$ is a diagonal matrix, it means that the integrals $\int_{W_n} f_i(X(t))\, dt,$ $i=1,\dots, s,$ are ``asymptotically independent''. Moreover, choosing the functions $f_1,\dots, f_s$ this way reduces the computational effort of the test from \cite{BST} and \cite[Sec.\ 3.5.2]{Kar} and we hope that this also improves the statistical properties of the test. 

First let us remark that the Gram-Schmidt orthogonalisation procedure is in principle suitable for this purpose. Consider a space of functions $V$ such that
\begin{equation} \int_{\R{d}} \big|\Cov\big(f(X(0)), g(X(t)) \big)\big|\, dt <\infty\label{e:asyFin} \end{equation}
for all $f,g\in V$. Then we see from Lemma \ref{l:asyForms} below that
\[ \langle\cdot,\cdot\rangle: V \times V\to \mathbb{R}, \ (f,g) \mapsto  \int_{\R{d}} \Cov\big(f(X(0)), g(X(t)) \big)\, dt \]
 is a symmetric non-negative definite bilinear form. Start with vectors $g_1,\dots, g_s \in V$ that fulfill the strengthened linear independence condition: If $\langle \sum_{i=1}^s \lambda_ig_i, \sum_{i=1}^s \lambda_ig_i \rangle=0$ for some $\lambda_1,\dots, \lambda_s\in\mathbb{R}$ then $\lambda_1=\dots=\lambda_s=0$. Define inductively
\begin{align*}
f_1&:=g_1\\
f_j&:=  g_j - \sum_{k=1}^{j-1} \frac{\langle g_j, f_k\rangle }{\langle f_k,f_k\rangle} f_k, \quad j=2, \dots, s.
\end{align*}
The resulting functions $f_i,\, i=1,\dots,s$ fulfill \eqref{e:Cov0}. 

 The strengthened linear independence condition is the reason why the Gram-Schmidt procedure works only ``in principle''. This condition is hard to check and in Example \ref{E:Sigma0} we construct a random field $(X(t))_{t\in\R{d}}$ such that there are no functions $g_1,\dots, g_s$ fulfilling it. 

\begin{lemma}\label{l:asyForms}
Let $(X(t))_{t\in\mathbb{R}^d}$ be a measurable, stationary random field and let $f,g:\mathbb{R}\to\mathbb{R}$ be two measurable functions fulfilling $\mathbb{E}\, f(X(0))^2<\infty$, $\mathbb{E}\, g(X(0))^2<\infty$ and \eqref{e:asyFin}. Let $(W_n)_{n\in\mathbb{N}}$ be a VH-growing sequence. Then
\[\lim_{n\to\infty} \frac{1}{\lambda_d(W_n)} \Cov\Big( \int_{W_n} f(X(t))\, dt, \int_{W_n} g(X(t))\, dt  \Big) =\langle f,g \rangle.\] 
\end{lemma} 
\prf We have 
\begin{align*}
\lim_{n\to\infty} \frac{1}{\lambda_d(W_n)}& \Cov\Big( \int_{W_n} f(X(t))\, dt, \int_{W_n} g(X(t))\, dt  \Big)\\
 & = \lim_{n\to\infty} \frac{1}{\lambda_d(W_n)}  \int_{W_n}\int_{W_n}\Cov\Big( f(X(s)),  g(X(t)) \Big)\, dt \, ds\\
 & = \lim_{n\to\infty} \frac{1}{\lambda_d(W_n)}  \int_{W_n}\int_{W_n-s}\Cov\Big( f(X(s)),  g(X(t+s)) \Big)\, dt\, ds \\
 & = \lim_{n\to\infty} \frac{1}{\lambda_d(W_n)}  \int_{W_n}\int_{\R{d}} \mathbf{1}_{W_n}(s+t)\Cov\Big( f(X(0)),  g(X(t)) \Big)\, dt\, ds \\
 & = \lim_{n\to\infty} \frac{1}{\lambda_d(W_n)}  \int_{\R{d}}\Cov\Big( f(X(0)),  g(X(t)) \Big)\int_{W_n} \mathbf{1}_{W_n}(s+t)\, ds \, dt\\
 & = \lim_{n\to\infty}  \int_{\R{d}}\Cov\Big( f(X(0)),  g(X(t)) \Big)\frac{\lambda_d(W_n\cap (W_n-t))}{\lambda_d(W_n)}\, dt.
\end{align*}
Now the assumption that $(W_n)_{n\in\mathbb{N}}$ is VH-growing implies 
\[ \lim_{n\to\infty} \frac{\lambda_d(W_n\cap (W_n-t))}{\lambda_d(W_n)} =1;\]
see \cite[Lemma 1.2, p.\ 172]{BS}.
By assumption \eqref{e:asyFin} we can apply the Dominated Convergence Theorem and thus the above expression equals $\langle f,g \rangle$. \qed\smallskip

Now we want to study further examples of functions $f_1,\dots, f_s$ which make the asymptotic covariance matrix $\Sigma$ a diagonal matrix. For a certain class of random fields which are defined using the L\'evy-Meixner systems introduced in Section \ref{ss:Meix} we get quite explicit examples.

\begin{exa}\label{e:Meix}
Let $\Psi_\lambda,\, \lambda>0$, be a L\'evy-Meixner system. Let $\Lambda$ be a random signed measure such that 
\begin{enumerate}[(i)]
\item $\Lambda(B)\sim \Psi_{\lambda_d(B)}$ for every Borel set $B\subseteq \mathbb{R}^d$ with $0< \lambda_d(B) <\infty$ and
\item $\Lambda(B_1)$ and $\Lambda(B_2)$ are independent, if $B_1$ and $B_2$ are disjoint.  
\end{enumerate}

Such a random measure exists due to the L\'evy noise construction. 

Define a random field $X$ by
$ X(t) :=  \Lambda(B+t),\ t\in\mathbb{R}^d,$ for a fixed Borel set $B\subseteq \mathbb{R}^d$ with $0< \lambda_d(B) <\infty$. 
Let $(Q_n(\cdot;\lambda))_{n\in\mathbb{N}_0}$ be the system of orthogonal polynomials w.r.t.\ $\Psi_\lambda$ appearing in Definition \ref{d:Meixner}. Then Theorem \ref{T:CovMeix}b) implies 
\[ \Cov\big( Q_n(X(t_1);\lambda_d(B)), Q_m(X(t_2); \lambda_d(B)) \big) =0\]
if $n\ne m$ and $t_1,t_2\in \mathbb{R}^d$.   
\end{exa}


\subsection{Non-degenerateness}

Now we would like to examine whether the asymptotic variance occurring in Theorem \ref{T:multCLTmix} is positive (more precisely: whether all diagonal entries of the asymptotic covariance matrix $\Sigma$ are positive). If this asymptotic variance is zero, we have not chosen the optimal normalization constant in Theorem \ref{T:multCLTmix}. Moreover, for the application of the Gram-Schmidt procedure in Section \ref{ss:Diag} it is important to know for which functions $f$ the asymptotic variance vanishes.

\begin{prop}\label{e:Sigma_reg}
 Let under the assumptions of Theorem \ref{T:multCLTmix} the field $(X(t))_{t\in\mathbb{R}^d}$ be either
\begin{itemize}
\item centered Gaussian with non-negative covariance function or
\item one of the fields constructed in Example \ref{e:Meix}.
\end{itemize}
Then we have:
\begin{enumerate}[1)]
\item If $\Var(\sum_{i=1}^s u_if_i(X(0)))>0$ for some $u=(u_1,\dots, u_s)\in\R{s}$ then $u^T\Sigma u>0$.
\item If $\Var(f_i(X(0)))>0$ then the $i$-th diagonal entry of $\Sigma$ is positive.
\item If $\Var\big(\sum_{i=1}^s u_if_i(X(0))\big)=0$ implies $u_1=\dots=u_s=0$ then $\Sigma$ is positive definite. 
\end{enumerate}
\end{prop}

\prf We only have to show 1), since 2) and 3) are immediate consequences. 

Let $\Psi_\lambda$ be the distribution of $X(0)$. The system $(Q_n(x;\lambda))_{n\in\mathbb{N}_0}$ of orthogonal polynomials of $\Psi_\lambda$ is complete (see Satz 4.2 and Satz 5.2 in \cite[\S\ II]{Fr}). Hence there are constants $c_n\in\mathbb{R}, n\in\mathbb{N}_0,$ with
\[ \sum_{i=1}^s u_if_i(x)=\sum_{n=0}^\infty c_nQ_n(x;\lambda), \]
where the series is $\mathcal{L}^2(\Psi_\lambda)$-convergent. By Theorem \ref{T:CovMeix}b) we have
\[\Cov\big(Q_n(X(0);\lambda), Q_m(X(t);\lambda)\big) =0, \quad t\in\mathbb{R}^d,\ m\ne n.\]
So
\begin{align*}
u^T \Sigma u &= \int_{\R{d}} \Cov\Big(\sum_{i=1}^s u_i f_i(X(0)),\sum_{i=1}^s u_i f_i(X(t))\Big) \, dt\\
&= \int_{\R{d}}\sum_{n=0}^\infty \sum_{m=0}^\infty c_nc_m \Cov\big( Q_n(X(0);\lambda), Q_m(X(t);\lambda) \big) \, dt\\
& =\int_{\R{d}}\sum_{n=0}^\infty c_n^2 \Cov\big( Q_n(X(0);\lambda), Q_n(X(t);\lambda) \big) \, dt. 
\end{align*}
For fixed $t\in\mathbb{R}^d$ and all $n\in\mathbb{N}_0$ we have 
\[ \Cov\big( Q_n(X(0);\lambda), Q_n(X(t);\lambda) \big) \ge 0.\]
Indeed, there are independent random variables $Y_1, Y_3 \sim \Psi_{\lambda_1}$ and $Y_2\sim \Psi_{\lambda_2}$ for appropriate constants $\lambda_1,\lambda_2\in(0,\infty)$ with $\lambda_1+\lambda_2=\lambda$ such that $(X(0), X(t))=(Y_1+Y_2,Y_2+Y_3)$ in distribution. In case of a Gaussian random field choose $Y_2\sim\mathcal{N}(0, \Cov(X(0),X(t)) )$ and $Y_1,Y_3 \sim \mathcal{N}(0, \Var(X(0)) - \Cov(X(0),X(t)) )$. In the case of a random field constructed in Example \ref{e:Meix}, put $Y_1:= \Lambda(B \setminus (B+t))$,  $Y_2 :={\Lambda((B+t)\cap B)}$, and $Y_3:={\Lambda((B+t)\setminus B)}$, where $\Lambda$ is the random measure from Example \ref{e:Meix}. Let $(Q_n(x;\lambda_i))_{n\in\mathbb{N}_0}$ be the system of orthogonal polynomials of $\Psi_{\lambda_i}$, $i=1,2$. From Theorem \ref{T:CovMeix}a) we get
\begin{align*}
\Cov\big( Q_n(X(0);\lambda),& Q_n(X(t);\lambda) \big)\\
 &= \Cov\Big( \sum_{i=0}^n\binom{n}{i} Q_{n-i}(Y_1;\lambda_1)Q_i(Y_2;\lambda_2), \sum_{j=0}^n\binom{n}{j} Q_{n-j}(Y_3;\lambda_1)Q_j(Y_2;\lambda_2) \Big)\\
&= \sum_{i=0}^n\binom{n}{i} \sum_{j=0}^n\binom{n}{j} \mathbb{E}[Q_{n-i}(Y_1;\lambda_1)]\cdot \mathbb{E}[Q_{n-j}(Y_3;\lambda_1)] \cdot \Cov( Q_i(Y_2;\lambda_2), Q_j(Y_2;\lambda_2))\\
& = \Cov( Q_n(Y_2;\lambda_2), Q_n(Y_2;\lambda_2))\\
&\ge 0. 
\end{align*}
Moreover, for $n>0$ we have $\Var(Q_n(X(0);\lambda))>0$. Since $(Q_n(X(t)))_{t\in\mathbb{R}^d}$ is a measurable, stationary and square-integrable random field, one can show -- using arguments from the proof of \cite[Prop.\ 3.1]{Roy} -- that it is continuous in $2$-mean. So for $t$ sufficiently close to $0$ we get
\[\Cov\big( Q_n(X(0);\lambda), Q_n(X(t);\lambda) \big) > 0. \]
Since $c_n=0$ for all $n>0$ would imply $\sum_{i=1}^s u_i f_i(x)=c_0$ for $\Psi_\lambda$-a.a.\ $x\in\mathbb{R}$, which contradicts the assumption, we have $c_{n_0}\ne 0$ for at least one $n_0>0$ and hence
\begin{align*}
\int_{\R{d}}\sum_{n=0}^\infty c_n^2 \Cov\big( Q_n(X(0);\lambda), Q_n(X(t);\lambda) \big) \, dt 
				&= \sum_{n=0}^\infty c_n^2\int_{\R{d}} \Cov\big( Q_n(X(0);\lambda), Q_n(X(t);\lambda) \big) \, dt \\
				& \ge c_{n_0}^2\int_{\R{d}} \Cov\big( Q_{n_0}(X(0);\lambda), Q_{n_0}(X(t);\lambda) \big) \, dt  >0. 
\end{align*}
So $u^T\Sigma u>0$.\qed\smallskip


However, the asymptotic variance is not always positive. Indeed, we can construct a non-trivial field fulfilling the assumptions of Theorem \ref{T:multCLTmix}, such that the asymptotic variance vanishes for all functions $f$ simultaneously.

\begin{exa}\label{E:Sigma0}
Let $Y=(\xi_i)_{i\in\mathbb{N}}$ be a stationary Poisson process on $\R{d}$ with intensity $1$. Consider the \emph{Voronoi mosaic} generated by $Y$, i.e.\ the system of sets
\[ C(\xi_i,Y):=\{ v\in\R{d}\mid \|v-\xi_i\| \le \| v-\xi_j\|, j\ne i \}. \]
Now every point $t\in\R{d}$ is contained in exactly one of the cells $C(\xi_i,Y)$ with probability one; we denote the corresponding point of the point process by $\xi_{i_t}$.  Now we define 
\[ X(t):=\lambda_d(\{v\in C(\xi_{i_t},Y) \mid \|v-\xi_{i_t}\| \le \|t-\xi_{i_t}\|\} )/\lambda_d(C(\xi_{i_t},Y)),\quad t\in\mathbb{R}^d. \]

\begin{figure}
\begin{center}
\begin{tikzpicture}
\filldraw[color=lightgray] (0,0) circle(1.4cm);
\filldraw[color=white] (-2,0) -- (1, -2) -- (-2, -2);
\draw (-2,2) -- (2,2) -- (1,-2) -- (-2,0) -- (-2,2) -- (-2.3, 2.3);
\draw (-2,0) -- (-2.2, -0.5);
\draw (2,2) -- (2.3, 2.3); 
\draw (1,-2) -- (1.1,-2.3); 
\filldraw (0,0) circle(2pt);
\draw (0.4 ,0) node {$\xi_{i_t}$};
\filldraw (-1,-2) circle(2pt);
\draw (-0.6 ,-2) node {$\xi_{j}$};
\filldraw (-1,1) circle(1.5pt);
\draw (-0.8, 1) node {$t$}; 
\end{tikzpicture}
\end{center}

The value of the field $X$ in $t$ is the fraction of the volume of the Voronoi cell covered by the grey circle. 
\caption{The definition of $X(t)$ in Example \ref{E:Sigma0}}\label{F:Xdef}
\end{figure}

So each point $t\in\mathbb{R}^d$ is assigned the proportion of the Voronoi cell it lies in that is closer to the nucleus than $t$, see Figure \ref{F:Xdef}. 

Now the field $(X(t))_{t\in\mathbb{R}^d}$ is obviously measurable and strictly stationary. Also other assumptions of Theorem \ref{T:multCLTmix} are fulfilled for measurable and bounded functions $f_1,\dots, f_s: [0,1] \to \mathbb{R}$. 
Indeed, assumption (i) is trivial. In order to check assumption (iii), we let $I,J\subseteq \R{d}$ be two sets lying in halfspaces separated by a strip of width $r\ge 15$ and fulfilling $V_d((I\cup J)+B^d)\le \gamma$. Then $I\cup J$ can intersect at most $N:=\lfloor\gamma/\sqrt{d}^d\rfloor$ cells of the lattice $\frac{1}{\sqrt{d}} \mathbb{Z}^d$ and such a cell can be covered by one ball of radius $R:= r/30$. Hence there are $N$ balls of radius $R$ covering $I\cup J$. Say $I\subseteq \bigcup_{i=1}^M B_R(x_i)$ and $J\subseteq \bigcup_{i=M+1}^N B_R(x_i)$. There are configurations of $Y_{|B_{15R}(x_i)}$ that fully determine the field $(X(t))_{t\in B_R(x_i)}$ and there are configurations for which $(X(t))_{t\in B_R(x_i)}$ also depends on $Y_{|\mathbb{R}^d\setminus B_{15R}(x_i)}$, where $X_{|M}$ and $f_{|M}$ denote the restriction of a stochastic process $X$ and a function $f$ to a subset $M$ of its index set or its domain. In order to make this precise, we have to observe that every locally finite counting measure $y$ on $\R{d}$ induces a function $x:\R{d}\to \R{}$ by the construction explained in the beginning of this example. So let $E_{R,i}$ denote the event that there is a point configuration $y$ on $\mathbb{R}^d$ with $y_{|B_{15R}(x_i)}=Y_{|B_{15R}(x_i)}$ for which the induced function $x$ does not coincide with $X$ on $B_R(x_i)$, $x_{|B_{R}(x_i)}\ne X_{|B_{R}(x_i)}$. It is shown in \cite[p.\ 515]{SW} that 
\[ \mathbb{P}(E_{R,i}) \le (1+m_d)e^{-\kappa_dR^d},\ i=1,\dots, N,\]
for some constant $m_d$ depending only on $d$. 
For any set $A\in\sigma(X_{\cup_{i=1}^MB_{R}(x_i)})$ and $B\in\sigma(X_{\cup_{i=M+1}^NB_{R}(x_i)})$ the events $A \cap \big(\bigcup_{i=1}^M E_{R,i}\big)^c$ and $B\cap \big(\bigcup_{i=M+1}^N E_{R,i}\big)^c$ are independent and therefore
\begin{align*}
 |\mathbb{P}(A&\cap B) - \mathbb{P}(A)\mathbb{P}(B)| \\
&\le   \bigg|\mathbb{P}\bigg(A\cap B \cap \Big(\bigcup_{i=1}^N E_{R,i}\Big)^c\bigg) - \mathbb{P}\bigg(A \cap \Big(\bigcup_{i=1}^M E_{R,i}\Big)^c\bigg) \cdot \mathbb{P}\bigg(B\cap \Big(\bigcup_{i=M+1}^N E_{R,i}\Big)^c\bigg)\bigg| 
+  2 N \mathbb{P}(E_{R,1}) 
\\
& \le 2N (1+m_d)e^{-\kappa_dR^d}.
\end{align*}
Hence 
\[\alpha_\gamma (r) \le 2 \tfrac{\gamma}{d^{d/2}} (1+m_d)e^{-\kappa_dr^d/30^d}\]
and thus condition (iii) is fulfilled. 

Moreover, we have
\begin{equation}
 \int_{\mathbb{R}^d} \big|\Cov\big(f(X(0)), g(X(t))\big)\big| \, dt <\infty
\label{e:ii} \end{equation}
for any bounded, measurable functions $f,g:[0,1]\to\mathbb{R}$. 
For any point $t\in\mathbb{R}^d$ let $E_{(0)}$ and $E_{(t)}$ be the events $E_{R,i}$ defined above with $x_i=0$ and $x_i=t$ respectively and with $R=\|t\|/30$. Put
\[ S:= \sup\{ |f(x)|, |g(x)| \mid x\in[0,1] \}. \]
We have 
\begin{align*}
 \big|\Cov\big(f(X(0))\mathbf{1}_{E_{(0)}}, g(X(t))\mathbf{1}_{(E_{(t)})^c}\big) \big|
&=\begin{aligned}[t] \big| \mathbb{E}\big[ f(X(0))\mathbf{1}_{E_{(0)}} &\cdot  g(X(t)) \mathbf{1}_{(E_{(t)})^c} \big]\\ & -  \mathbb{E}\big[ f(X(0))\mathbf{1}_{E_{(0)}}\big] \cdot\mathbb{E}\big[  g(X(t))\mathbf{1}_{(E_{(t)})^c} \big] \big| \end{aligned}
\\
 & \le  \mathbb{E}\big[S^2\cdot \mathbf{1}_{E_{(0)}}\big]  +  \mathbb{E}\big[S\cdot \mathbf{1}_{E_{(0)}}\big] \cdot S\\
&= 2S^2\cdot \mathbb{P}(E_{(0)}),
\end{align*}
the same way 
\begin{align*}
 \big| \Cov\big(f(X(0)) \mathbf{1}_{(E_{(0)})^c}, g(X(t))\mathbf{1}_{E_{(t)}}\big) \big| &\le 2S^2\cdot  \mathbb{P}( E_{(t)}) 
\end{align*}
and 
\begin{align*}
 \big| \Cov\big(f(X(0))\mathbf{1}_{E_{(0)}}, g(X(t))\mathbf{1}_{E_{(t)}}\big) \big|
&= \begin{aligned}[t] \big| \mathbb{E}\big[ f(X(0))\mathbf{1}_{E_{(0)}} &\cdot  g(X(t)) \mathbf{1}_{E_{(t)}} \big]\\ & -  \mathbb{E}\big[ f(X(0))\mathbf{1}_{E_{(0)}}\big] \cdot\mathbb{E}\big[  g(X(t))\mathbf{1}_{E_{(t)}} \big] \big| \end{aligned} \\
 & \le \mathbb{E}\big[ S^2 \cdot\mathbf{1}_{E_{(0)}} \cdot\mathbf{1}_{E_{(t)}}\big]  +  \mathbb{E}\big[ S\cdot \mathbf{1}_{E_{(0)}}\big] \cdot \mathbb{E}\big[  S\mathbf{1}_{E_{(t)}} \big]\\
&= 2S^2 \cdot \mathbb{P}(E_{(0)})\cdot\mathbb{P}( E_{(t)}). 
\end{align*}
Since, moreover, the random variables $f(X(0))\mathbf{1}_{(E_{(0)})^c}$ and $g(X(t))\mathbf{1}_{ ( E_{(t)})^c}$ are independent, we get
\begin{align*}
 \big|\Cov\big(f(X(0)), g(X(t))\big)\big| &\le \big|\Cov\big(f(X(0))\mathbf{1}_{(E_{(0)})^c}, g(X(t))\mathbf{1}_{(E_{(t)})^c}\big)\big| + 6S^2 \cdot \mathbb{P}(E_{(0)} )\\ 
&\le 6S^2 (1+m_d)e^{-\kappa_d\|t\|^d/30^d}. 
\end{align*}
Hence \eqref{e:ii} holds and thus all assumptions of Theorem \ref{T:multCLTmix} are fulfilled.

However, for any measurable, bounded function $f:\mathbb{R}\to\mathbb{R}$ the asymptotic variance is zero:
\[ \int_{\R{d}} \Cov\big(f(X(0)), f(X(t))\big)  \, dt =0. \]
 At first we derive bounds for
$\int_{[-r,r]^d} \Cov\big(f(X(0)), f(X(t)) \big) \, dt$
for any $r>0$. Put $x_1:=0$ and choose $N-1$ balls of radius $R:=r/30$ covering the boundary of $[-r,r]^d$, say $B_R(x_2), \dots, B_R(x_N)$. Clearly, $N$ can be chosen independent of $r$. Let $C_i$ denote the union of all cells that lie completely in $[-r,r]^d$ and $C_b$ the union of all cells intersecting the boundary of $[-r,r]^d$. Then
\begin{align*} 
\Big| \int_{[-r,r]^d}\Cov\big( f(X(0)), f(X(t)) \big)\, dt \Big|
 \le &
	\Big| \Cov\Big( f(X(0)) \cdot \mathbf{1}_{( E_{R,1})^c}, \int_{C_i} f(X(t)) \, dt  \cdot \mathbf{1}_{(\cup_{i=2}^N E_{R,i})^c}\Big) \\
	&\quad+ \Cov\Big( f(X(0)) \cdot \mathbf{1}_{( E_{R,1})^c}, \int_{C_b\cap [-r,r]^d} f(X(t)) \, dt  \cdot \mathbf{1}_{(\cup_{i=2}^N E_{R,i})^c}\Big) \Big| \\
	+& 2S^2\cdot (2r)^d N\mathbb{P}(E_{R,1}).
\end{align*}
Now $f(X(0)) \cdot \mathbf{1}_{(E_{R,1})^c}$ and $\int_{C_b\cap [-r,r]^d} f(X(t)) \, dt  \cdot \mathbf{1}_{(\cup_{i=2}^N E_{R,i})^c}$ are independent.

In order to compute $\int_{C_i} f(X(t)) \, dt$, we calculate $\int_{C(\xi_{i},Y)} f(X(t))\, dt$ for every cell $C(\xi_{i},Y)$ of the mosaic by algebraic induction. First let $f=\mathbf{1}_{[0,b]}, b\in[0,1]$. We introduce the function
\[ g: [0,\infty)\to\mathbb{R},\, s\mapsto \lambda_d(\{v\in C(\xi_i,Y)  \mid \|v-\xi_i\| \le s\}). \]  
By \cite[Lemma 5, Lemma 2(i) and Lemma 1]{Sta}, it is continuous and strictly monotonically increasing on $\{s\in [0,\infty) \mid g(s)<\lambda_d(C(\xi_i,Y)) \}$. Hence the restriction of $g$ to $[0,s^*]$, where
\[s^*:=\min\{ s\ge 0 \mid g(s)=\lambda_d(C(\xi_i,Y))\}, \]
is  bijective. We have
\[ X(t) = g(\|t-\xi_i\|)/\lambda_d(C(\xi_i,Y)),\quad t\in C( \xi_i,Y), \]
and therefore
\begin{align*}
 \int_{C(\xi_{i},Y)} \mathbf{1}_{[0,b]}(X(t))\, dt &= \lambda_d(\{t\in C(\xi_i,Y) \mid g(\|t-\xi_i\|) \le b\cdot \lambda_d(C(\xi_i,Y)) \}) \\
&= \lambda_d(\{t\in C(\xi_i,Y) \mid \|t-\xi_i\| \le g^{-1}(b\cdot \lambda_d(C(\xi_i,Y))) \}) \\
&= g( g^{-1}(b\cdot \lambda_d(C(\xi_i,Y)))) \\
&= b\cdot \lambda_d(C(\xi_i,Y)). 
\end{align*}
The system
\[ \Big\{ B\in\mathcal{B}([0,1]) \Big|  \int_{C(\xi_{i},Y)} \mathbf{1}_{B}(X(t))\, dt =  \lambda(B)\cdot \lambda_d(C(\xi_i,Y)) \Big\} \]
is clearly a Dynkin system and hence coincides with Borel-$\sigma$-algebra $\mathcal{B}([0,1])$. Now it is easy to see that the equality
\[ \int_{C(\xi_{i},Y)} f(X(t))\, dt =  \int_{[0,1]} f(x)\, dx \cdot\lambda_d(C(\xi_i,Y)) \]
extends to linear combinations of indicator functions, thus to non-negative measurable functions $f$ and finally to integrable functions $f:[0,1]\to\mathbb{R}$. 

So
\begin{align*}
 	\Cov\Big( f(X(0)) \cdot \mathbf{1}_{(E_{R,1})^c}, \int_{C_i} f(X(t)) \, dt  &\cdot \mathbf{1}_{(\cup_{i=2}^N E_{R,i})^c}\Big) \\
	&= 	\Cov\Big( f(X(0)) \cdot \mathbf{1}_{(E_{R,1})^c}, \int_{[0,1]} f(x)\, dx \cdot\lambda_d(C_i)  \cdot \mathbf{1}_{(\cup_{i=2}^N E_{R,i})^c}\Big).
	\end{align*}
However, if $\bigcup_{i=2}^N E_{R,i}$ does not hold, then the boundary of $C_i$ and hence $\lambda_d(C_i)$ is completely determined by $Y_{|\cup_{i=2}^N B_{15R}(x_i)}$. So $f(X(0)) \cdot \mathbf{1}_{(E_{R,1})^c}$ and $\lambda_d(C_i)  \cdot \mathbf{1}_{(\cup_{i=2}^N E_{R,i})^c}$ are independent and therefore we conclude 
\[
\Big| \int_{[-r,r]^d}\Cov\big( f(X(0)), f(X(t)) \big)\, dt\Big|  \le 2S^2\cdot 2^d r^dN\mathbb{P}(E_{R,1})
 \le 2^{d+1}S^2r^dN(1+m_d)e^{-\kappa_dr^d/30^d}. 
\]
	Thus
	\[ \int_{\R{d}}\Cov\big( f(X(0)), f(X(t)) \big)\, dt =0. \] 
\end{exa}

\section{Functional central limit theorem}\label{s:func}

In this section we are going to extend Theorem \ref{T:multCLTmix} to a functional central limit theorem.

We let $V$ denote the vector space of all Lipschitz continuous functions $\mathbb{R} \to \mathbb{R}$ and equip $V$ with the norm defined by 
\[ \|f\|:= \Lip f + |f(0)|,\]
where 
\[ \Lip f:= \sup_{x,y\in\mathbb{R},  x\ne y} \frac{|f(x)-f(y)|}{|x-y|} \]
is the Lipschitz constant of $f$.

For a measurable random field $(X(t))_{t\in\mathbb{R}^d}$ and a sequence $(W_n)_{n\in\mathbb{N}}$ of compact subsets of $\mathbb{R}^d$ we define random functionals by
\[\Phi_n: V\to \mathbb{R}, f \mapsto \frac{\int_{W_n} f(X(t)) \, dt - \lambda_d(W_n) \mathbb{E}[f(X(0))] }{\sqrt{\lambda_d(W_n)}}. \]
Under weak regularity conditions, these functionals take values in the dual space $V^*$ of $V$.
\begin{lemma} 
If $\int_{W_n} |X(t)|\, dt<\infty$ a.s.\ for all $n\in\mathbb{N}$ and $\mathbb{E}[|X(0)|]<\infty$, then the paths of $\Phi_n,\, n\in \mathbb{N},$ are linear and continuous. 
\end{lemma} 
\prf The linearity is readily checked. In order to prove the continuity, it suffices to show that the paths are bounded on the set $\{ f\in V \mid \|f\|\le 1 \}$. This holds since $|f(x)| \le \|f\| \cdot (|x|+1)$ for all $x\in\mathbb{R}$ and therefore
\[ |\Phi_n(f)| \le \|f\| \cdot \frac{\int_{W_n} (|X(t)|+1) \, dt + \lambda_d(W_n) \mathbb{E}[|X(0)|+1] }{\sqrt{\lambda_d(W_n)}} .\qed \]
We consider $V^*$ with the weak topology, i.e.\ the topology of pointwise convergence. As usual, we say that a sequence $(\Phi_n)_{n\in\mathbb{N}}$ of $V^*$-valued random variables converges in distribution to a $V^*$-valued random variable $\Phi$ if
\[ \lim_{n\to\infty} \mathbb{E}\, h(\Phi_n) = \mathbb{E}\, h(\Phi) \]
for all bounded continuous functions $h:V^* \to \mathbb{R}$. 

\begin{theorem}\label{t:FCLT}
Let $(X(t))_{t\in \R{d}}$ be a stationary and measurable random field. Assume that 
 there are constants $n\in\mathbb{N}$, $\delta>4$ and $C, l>0$ with $n/d>\max\{l + \delta/(\delta-2),\delta/(\delta-4),3\}$ such that $\alpha_\gamma(r)\le C r^{-n}\gamma^l$ for all $\gamma\ge 2\kappa_d,\, r>0$ and $\mathbb{E} |X(0)|^\delta<\infty$. 
Let $(W_n)_{n\in\mathbb{N}}$ be a VH-growing sequence of compact subsets of $\mathbb{R}^d$. Then $(\Phi_n)_{n\in\mathbb{N}}$ converges weakly in distribution to a centered Gaussian process $\Phi$ on $V$ with covariance
\begin{equation} \Cov( \Phi(f), \Phi(g)) =\int_{\R{d}} \Cov( f(X(0)), g(X(t)) )\, dt.\label{e:CovPhi} \end{equation}
\end{theorem}

The proof of this theorem relies on the following result, which is a special case of Satz 5 of \cite{Op}:  
\begin{theorem}
Let $V$ be a locally convex, metrizable space. A sequence of $V^*$-valued random elements $(\Phi_n)_{n\in\mathbb{N}}$ converges to a $V^*$-valued random element $\Phi$ iff for all finite collections $\{f_1,\dots, f_m\}$ the sequence $(\Phi_n(f_1), \dots, \Phi_n(f_m))_{n\in\mathbb{N}}$ converges in distribution to $(\Phi(f_1),\dots, \Phi(f_m))$.
\end{theorem}
So, in order to prove Theorem \ref{t:FCLT}, we have to prove the convergence of the finite-dimensional random vectors and we have to show that there is a centered Gaussian process $\Phi$ with covariance \eqref{e:CovPhi} and that this process has a version with linear and continuous paths.

The first point is simply checking the conditions of Theorem \ref{T:multCLTmix}. All assumption of Theorem \ref{T:multCLTmix} except assumption (ii) are directly assumed in Theorem \ref{t:FCLT}.  
By \cite[p.\ 9, Theorem 3(1)]{Dou} we have 
\begin{align*}  |\Cov(f_i(X(0)),f_j(X(t)))| &\le 8\alpha_{2\kappa_d}(\|t\|_\infty)^{1/3} \mathbb{E}[f_i(X(0))^3]^{1/3} \mathbb{E}[f_j(X(t))^3]^{1/3}\\
& \le  C'\|t\|_\infty^{-d-\epsilon} \mathbb{E}[f_i(X(0))^3]^{1/3} \mathbb{E}[f_j(X(0))^3]^{1/3}
\end{align*}
for any Lipschitz continuous functions $f_i,f_j:\mathbb{R}\to \mathbb{R}$ and some constants $C', \epsilon>0$. Since $\mathbb{E}[f_i(X(0))^3]<\infty$ by Lemma \ref{l:LipE}(i) below, 
we have
\[\int_{\mathbb{R}^d} |\Cov(f_i(X(0)),f_j(X(t)))| \, dt <\infty,\]
which is assumption (ii) of Theorem \ref{T:multCLTmix}.
Thus the assumptions of Theorem \ref{T:multCLTmix} are fulfilled and therefore the finite-dimensional distributions converge appropriately.

 




In order to show the existence of the process $\Phi$, we employ the theory of $GB$- and $GC$-sets (see e.g. \cite[Chapter 2]{Dud99}). The \emph{isonormal process} on a Hilbert space $H$ with scalar product $\langle \cdot , \cdot \rangle$ is the centered Gaussian process $\Psi$ with covariance function 
\[\Cov(\Psi(f), \Psi(g)) =\langle f, g\rangle, \quad f,g\in H.\]
For a set $B\subseteq H$ and $\epsilon>0$ we define  $N(\epsilon)$ to be the minimal number $n$ of elements $f_1,\dots, f_n\in B$ such that for any element $f\in B$ there is an index $i\in \{1,\dots, n\}$ with $\sqrt{\langle f-f_i, f-f_i\rangle} < \epsilon$.
We obtain combining \cite[Theorem 2.6.1, Theorem 2.5.5(g)]{Dud99} and the remark above \cite[Lemma 2.5.3]{Dud99}:
\begin{lemma}\label{L:Psi_version} 
Let $H$ be a Hilbert space and $B\subseteq H$. If 
\begin{equation}
 \int_0^\infty (\log N(\epsilon))^{1/2}\, d\epsilon<\infty, 
\label{e:GCint} \end{equation}
then the isonormal process $\Psi$ has a version with paths that are bounded on $B$ and linear on $H$. 
\end{lemma}

We define a symmetric non-negative definite bilinear form $\langle \cdot, \cdot \rangle$ on $V$ by
\[ \langle f, g\rangle = \int_{\mathbb{R}^d} \Cov( f(X(0)), g(X(t)) )\, dt. \]
Now $\tilde H$, defined as the quotient space $V$ modulo the functions $f\in V$ with $\langle f, f \rangle =0$, is a pre-Hilbert space. Its completion $H$ is a Hilbert space, see \cite[Theorem 5.4.11]{DudRAP}. It is not clear whether $\tilde H=H$ holds. 


\begin{lemma}\label{L:GB}
For the canonical projection of the set
$ B:=\{f\in V\mid \| f\| \le 1\} $
onto $H$ inequality \eqref{e:GCint} holds. 
\end{lemma}

\begin{kor}
There is a centered Gaussian process $\Phi$ on $V$ with covariance function \eqref{e:CovPhi}. This process has a version whose paths are linear and continuous.
\end{kor}
\prf Let $\pi:V\to H$ be the canonical projection and injection map. According to Lemma \ref{L:GB} and Lemma \ref{L:Psi_version} the isonormal process on $H$ has a version $\Psi$ with paths which are linear on $H$ and bounded on $\pi(B)$. So $\Phi:=\Psi\circ \pi$ is a centered Gaussian process with covariance function \eqref{e:CovPhi} and linear paths. Since the paths are bounded on $B$, they are continuous.   \qed\smallskip

\prf[ of Lemma \ref{L:GB}]

We construct the functions $f_1, \dots, f_n$ like this: Fix $m\in \mathbb{N}$ and $c>0$. Then there are $n=2^{2m}$ Lipschitz-continuous functions $f$ with the following properties:
\begin{itemize}
\item $f(0)=0$
\item On each of the intervals $[(k-1)c, kc],\, k=-m+1,\dots, m,$ the function $f$ is either increasing with constant slope $1$ or decreasing with constant slope $-1$.
\item On the intervals $(-\infty,-mc]$ and $[mc, \infty)$, the function $f$ is constant.
\end{itemize}
Enumerate these functions as $f_1, \dots, f_n$. Now for each function $f\in B$ with $f(0)=0$ there is a function $f_i$ with $|f(kc)-f_i(kc)|\le c$ for each $k=-m,\dots, m$. We will prove this claim by induction for a fixed function $f$. Of course, $|f(kc)-f_i(kc)| =0 \le c $ for $k=0$ and all $i=1,\dots, n$. Assume that the existence of an $i\in\{1,\dots, n\}$ with $|f(kc)-f_i(kc)|\le c$ for $k=-M,\dots, M$ has been shown for some $M<m$. By the definition of the system $\{f_1,\dots, f_n\}$ it is clear that as soon as there is one function $f_i$ with this property, there are four such functions $f_{i_1},\dots, f_{i_4}$ such that on $[Mc,(M+1)c]$ the functions $f_{i_1}$ and $f_{i_3}$ are increasing and $f_{i_2}$ and $f_{i_4}$ are decreasing, while on $[(-(M+1)c, -Mc]$ the functions $f_{i_1}$ and $f_{i_2}$ are increasing and $f_{i_3}$ and $f_{i_4}$ are decreasing. Hence $|f((M+1)c)-f_i(Mc)|<2c$ for $i=i_1,i_2,$ and moreover $f_{i_1}((M+1)c)=f_{i_1}(Mc)+c$ and $f_{i_2}((M+1)c)=f_{i_2}(Mc)-c$ which implies $|f((M+1)c)-f_i((M+1)c)|\le c$ either for $i=i_1$ or for $i=i_2$. Of course, this holds also for $i=i_3$ resp.\ $i=i_4$. Moreover, in the same way we get $|f((-M-1)c)-f_i((-M-1)c)|\le c$ either for $i=i_1, i_2$ or for $i=i_3, i_4$. Since there is one index $i$ for which both inequalities hold, the induction step is completed.


Now 
choose some $r<\delta/2$ with $n/d>r/(r-2)$, which is equivalent to $n(1-2/r)>d$. Put $\gamma:=\min\{r, \delta-2r\}$, $c:=\epsilon^{1+\gamma/(4r)}$ and $m:=\lceil \tfrac{1}{\epsilon^{2-\gamma/(4r)}}\rceil$, where $\epsilon>0$ is sufficiently small. Then for each $f\in B$ there is an index $i\in \{1,\dots, n\}$ with $\langle f-f_i,f-f_i\rangle < \epsilon^2$. Indeed, since the bilinear form vanishes for constant functions, we may assume $f(0)=0$. Now let $i$ be the same index as chosen right above. We abbreviate $g:=f-f_i$. Then \cite[p.\ 9, Theorem 3(1)]{Dou} implies 
\[
\langle f-f_i, f-f_i\rangle = \int_{\mathbb{R}^d} \Cov\big(g(X(0)), g(X(t)) \big) \, dt   \le 8 (\mathbb{E}[|g(X(0))|^r])^{2/r} \cdot \int_{\R{d}} \alpha_{2\kappa_d}(\|t\|)^{(1-2/r)}\, dt.
\]
The integral is finite, since $\alpha_{2\kappa_d}(\|t\|)^{(1-2/r)} \le C(2\kappa_d)^l \|t\|^{-n(1-2/r)}$. Since $|g(x)|\le 2|x|$ for all $x\in\mathbb{R}$, Lemma \ref{l:LipE}(ii) below implies 
\begin{align*}
& \mathbb{E}[ |g(X(0))|^r \mathbf{1}_{\{X(0)\le -y\}}]\le \mathbb{E}[ 2^r|X(0)|^r \mathbf{1}_{\{X(0)\le -y\}}]\le D\cdot y^{-(\delta-r)}\le D\cdot y^{-(r+\gamma)} \quad\mbox{and}\\
& \mathbb{E}[ |g(X(0))|^r \mathbf{1}_{\{X(0)\ge y\}}] \le \mathbb{E}[ 2^r|X(0)|^r \mathbf{1}_{\{X(0)\ge y\}}] \le D\cdot y^{-(\delta-r)} \le D\cdot y^{-(r+\gamma)} 
\end{align*}
for $D:=2^r\mathbb{E}[X(0)^\delta]$ and $y\ge 1$. Moreover, $|g(kc)|<c$ for $k=-m,\dots, m$ and the fact that $g$ is Lipschitz continuous with Lipschitz constant at most $2$ imply $|g(x)|<2c$ for $x\in [-mc, mc]$. Hence
\begin{align*}
 \mathbb{E}[|g(X(0))|^r] 
&=  \mathbb{E}[|g(X(0))|^r\mathbf{1}_{\{X(0)\le -mc\}}] + \mathbb{E}[|g(X(0))|^r\mathbf{1}_{\{-mc < X(0)< mc\}}] + \mathbb{E}[|g(X(0))|^r\mathbf{1}_{\{X(0)\ge mc\}}] 
\\
&\le D\cdot (mc)^{-(r+\gamma)} + (2c)^r + D\cdot(mc)^{-(r+\gamma)}\\
&= 2D\cdot(\lceil 1/\epsilon^{2-\gamma/(4r)} \rceil \cdot \epsilon^{1+\gamma/(4r)})^{-(r+\gamma)} + 2^r\epsilon^{r+\gamma/4}\\
&\le 2D\cdot  \epsilon^{(1-\gamma/(2r))(r+\gamma)} + 2^r\epsilon^{r+\gamma/4}\\
&=: h(\epsilon).\end{align*}
Since $(1-\gamma/(2r))(r+\gamma)= r+\gamma/2-\gamma^2/(2r)>r$, we get $\lim_{\epsilon\to 0} h(\epsilon) /\epsilon^r = 0$ and so $\langle f-f_i, f-f_i\rangle$ is less than $\epsilon^2$ if $\epsilon>0$ is small enough. 

Thus we have shown $N(\epsilon)\le 2^{2m} \le 2^{2/\epsilon^{2-\gamma/(4r)}+2}$ for all $\epsilon\in(0,E_1)$ for some constant $E_1>0$. Therefore 
\begin{align*}
\int_0^{E_1} (\log N(\epsilon))^{1/2}\, d\epsilon &< \int_0^{E_1} (\log 2)^{1/2}\cdot (2/\epsilon^{2-\gamma/(4r)}+2)^{1/2} \, d\epsilon\\
& \le (\log 2)^{1/2} \cdot \int_0^{E_1} 2\cdot\max\{\epsilon^{-(1-\gamma/(8r))},1\} \, d\epsilon
<\infty.
\end{align*}
Since $\int_{E_1}^\infty (\log N(\epsilon))^{1/2}\, d\epsilon<\infty$ by \cite[p.\ 52]{Dud99}, this completes the proof of the lemma. \qed

\begin{lemma}\label{l:LipE}
Let $X$ be a random variable with $\mathbb{E}[|X|^r]<\infty$, $r\in\mathbb{R}$. 
\begin{enumerate}[(i)]
\item For every Lipschitz continuous function $f:\mathbb{R} \to \mathbb{R}$ we have $\mathbb{E}[|f(X)|^r]<\infty$.
\item For any $s<r$ and all $x>0$ we have
\[ \mathbb{E}[X^s\mathbf{1}_{\{X\ge x\}}] \le \mathbb{E}[|X|^r] \cdot x^{-(r-s)}. \]
\end{enumerate}
\end{lemma}
\prf (i) There are constants $a,b\in\mathbb{R}$ with $|f(x)|\le \max\{a|x|, b\}$. Hence
\[ \mathbb{E}[|f(X)|^r] \le \mathbb{E}[\max\{a^r|X|^r, b^r\}] \le \mathbb{E}[a^r|X|^r+ b^r] <\infty.\]
(ii) We have $X^s\mathbf{1}_{\{X\ge x\}} \le |X|^r \cdot x^{-(r-s)}$ for each fixed realisation. Taking expected values yields the assertion. \qed

\begin{rem} It remains an open question whether the mixing assumptions in Theorem \ref{t:FCLT} can be replaced by association assumptions. 
\end{rem}

\end{document}